\documentclass{article}

\usepackage[english]{babel}

\usepackage[utf8]{inputenc}
\usepackage[T1]{fontenc}

\usepackage{amssymb}
\usepackage{stmaryrd}
\usepackage{amsmath,amsfonts}
\usepackage[tikz]{bclogo}
\usepackage[top=4.34cm]{geometry}

\usepackage{lmodern}
\usepackage{textcomp}
\usepackage{mathtools, bm}
\usepackage{amssymb, bm}
\usepackage[unq]{unique}
\usepackage{enumerate}
\usepackage{bbm}

\usepackage[breaklinks]{hyperref}
\usepackage[numbers]{natbib}    

\usepackage{amsthm}
\usepackage{thmtools}
\usepackage{comment}
\usepackage{cleveref}

\newtheorem{theorem}{Theorem}[section]

\newtheorem{Question}[theorem]{Question}

\theoremstyle{definition}

\title{A note on Cops and Robbers, independence number, domination number and diameter}
\author{Jan Petr \and Julien Portier  \and Leo Versteegen\footnote{\{jp895, jp899, lvv23\}@cam.ac.uk, Centre for Mathematical Sciences, Wilberforce Road, Cambridge CB3 0WA, United Kingdom}}
\date{}

\begin{document}

\maketitle

\begin{abstract}
We study relations between diameter $D(G)$, domination number $\gamma(G)$, independence number $\alpha(G)$ and cop number $c(G)$ of a connected graph $G$, showing (i.) $c(G) \leq \alpha(G)-\lfloor \frac{D(G)-3}{2} \rfloor$, and (ii.) $c(G) \leq \gamma (G) - \frac{D(G)}{3} + O (\sqrt{D(G)})$.
\end{abstract}

\section{Introduction}

The game of Cops and Robbers, introduced by Nowakowski and Winkler \cite{Now} and Quillot \cite{Qui}, is played by two players on a finite undirected graph. The first player controls $k$ cops, the other one a single robber. At the beginning, the first player chooses a starting vertex for each of the cops. Afterwards, the second player chooses a starting vertex for the robber. From this point onwards, the players take alternating turns, beginning with the cops. During their turn, every cop moves to a vertex at distance at most $1$. In robber's turn, so does the robber. Multiple cops are allowed to stand on the same vertex during all stages of the game. The cops win if at any point the robber is on the same vertex as one of the cops. If on the other hand the robber evades the cops forever, the robber wins.

For a graph $G$, the minimal number $c(G)$ of cops such that there is a winning strategy for the cops is called the \emph{cop number of $G$}. Note that the cop number of a graph is equal to the sum of the cop numbers of its components. Therefore, we can restrict ourselves to connected graphs.

In this paper we investigate a question by Turcotte \cite{2k2free} about the relations between the cop number, independence number and domination number of a graph. It is easy to see that for any graph $G$, the cop number is bounded by the domination number $\gamma(G)$, i.e., the least size of a vertex set such that every other vertex in the graph has a neighbour in it, which in turn is at most the independence number $\alpha(G)$. Turcotte asked for which connected graphs we have $c(G)=\gamma(G)=\alpha(G)$. We show that this equality can only hold for graphs with diameter $D(G)$ at most $3$.

\begin{theorem}\label{AlphaMinus1}
Let $G$ be a connected graph of diameter $D(G) \geq 4$. Then $c(G) \leq \alpha(G)-1$.
\end{theorem}

It is also natural to ask for which connected graphs $c(G)=\gamma(G)$ holds. We show that this equality can only hold for graphs with diameter $D(G)$ at most $5$.

\begin{theorem}\label{GammaMinus1}
Let $G$ be a connected graph of diameter $D(G) \geq 6$. Then $c(G) \leq \gamma(G)-1$.
\end{theorem}

What is more, we show that the minimal difference between $\alpha(G)$ or $\gamma(G)$ on the one hand and $c(G)$ on the other grows at least linearly as a function of $D(G)$.

\begin{theorem}\label{prop:diam-idpd}
Let $G$ be a connected graph. Then $c(G) \leq \alpha(G)- \left \lfloor \frac{D(G)-3}{2} \right \rfloor$.
\end{theorem}

\begin{theorem}\label{prop:diam-domi}
Let $G$ be a connected graph. Then $c(G) \leq \gamma(G)-\dfrac{D(G)}{3} +O(\sqrt{D(G)})$.
\end{theorem}

The path $P_n$ on $n$ vertices shows that the above propositions are essentially best possible for connected graphs with large diameter as $D(P_n)=n-1$, $c(P_n)=1$, $\alpha(P_n)=\lceil \frac{n}{2}\rceil$ and $\gamma(P_n)=\lceil \frac{n}{3}\rceil$.

\section{Cop number and independence number}

For a vertex $u$ and a non-negative integer $i$, let the neighbourhood at distance $i$ of $u$ be $N_i(u)= \{ v \in V: \mathrm{dist}(u,v)=i \}$. By $N_{\geq i}(u)$ we will mean $\{v \in G: \mathrm{dist}(u,v) \geq i\}$. 
For a graph $G$, we say that a set $R \subset V(G)$ is \emph{dominated} by a positioning of cops if for every vertex $r$ of $R$, there is at least one cop in $N(r) \cup \{r\}$. \\
As we said in the introduction, for any graph $G$, $c(G) \leq \gamma(G) \leq \alpha(G)$ since putting cops on a dominating set of $G$ catches the robber for sure. We now prove \Cref{AlphaMinus1}, which can be seen as a strenghtening version of this inequality for graphs of diameter at least $4$.

\begin{proof}[Proof of \Cref{AlphaMinus1}]
Let $u$ and $v$ be vertices such that $\mathrm{dist}(u,v) = D(G)$, $N_i=N_i(u)$ and $N_{\geq i}=N_{\geq i}(u)$. Among all maximum independent sets of the subgraph
$G[N_{\geq 2}]$, choose $I$ such that $|I \cap N_3|$ is minimal.
Note that $I \cup \{u\}$ is a maximal independent set of $G$, hence $|I| \leq \alpha(G)-1$. We put a cop on every vertex of $I$. The robber must start on $\{ u \}$ or $N_1$ since $I$ dominates $N_{\geq 2}$. If $D(G) \geq 5$, then we are already done, because $I \cap N_{\geq 4}$ is non-empty, so we can take a cop from $N_{\geq 4}$ and bring it to $u$. This catches the robber since $N_2$ stays dominated by the cops through the whole process so the robber cannot leave $\{ u \} \cup N_1$.

Suppose therefore that $D(G)=4$. If $I \cap N_4$ is non-empty, then we can move a cop from it to $u$ as above, as the robber must stay in $\{ u \} \cup N_1$ since $N_2$ remains dominated by the cops. Otherwise, $I \cap N_4$ is empty, and therefore $I \cap N_3$ is non-empty. Take $w \in I \cap N_3$. We claim that
all the neighbours of $w$ in $N_2$ have a neighbour in $I\setminus \{w\}$.
Indeed, suppose there is a neighbour $w'$ of $w$ in $N_2$ that has no neighbour in $I'=I \setminus \{w\} \cup \{w'\}$. Then $I'$ is an independent set with $|I'|=|I|$ and $|I' \cap N_3|<|I \cap N_3|$, which contradicts the definition of $I$. We can therefore move the cop from $w$ to $u$, catching the robber as above.
\end{proof}

\begin{figure}[htbp]\centering
    			\includegraphics[height=3cm]{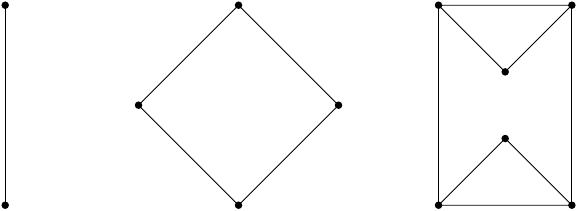}
    			\caption{Examples of connected graphs with $c(G)=\gamma(G)=\alpha(G)$ for $D(G) \in \{1,2,3\}$.}
\label{c=a}
	\end{figure}

This proposition is optimal in the sense that there are examples of connected graphs with $c(G)=\alpha(G)$ for diameter $D(G) \leq 3$, see \Cref{c=a}.

We can apply the same argument as in the proof of $\Cref{AlphaMinus1}$ to prove \Cref{prop:diam-idpd}.

\begin{proof}[Proof of \Cref{prop:diam-idpd}.]
Let $D(G)$ be the diameter of $G$. If $D(G) \leq 4$, the result is obvious. Otherwise, let $u$, $v$ be vertices such that $\mathrm{dist}(u,v)=D(G)$, and $N_i=N_i(u)$. For $0 \leq i \leq D(G)-2$, let $B_i$ be the induced subgraph of $G$ on $N_i \cup N_{i+1} \cup N_{i+2}$, $\alpha_i = \alpha(B_i)$, $j$ be such that $\alpha_j$ is maximal, and $S_j$ be a maximum independent set of $B_j$. We will show that $\alpha_j$ cops are enough to catch the robber. Note that, for $d_1=j$ and $d_2=D(G)-2-j$, we have: $ \left \lfloor \frac{d_1}{2} \right \rfloor + \left \lfloor \frac{d_2}{2} \right \rfloor + \alpha_j \leq \alpha(G)$ since $I=S_j \cup \{ a_0, a_1, \dots, a_{\lfloor d_1/2 \rfloor-1} \} \cup \{ b_0, b_1, \dots, b_{\lfloor d_2/2 \rfloor-1} \}$ is an independent set of $G$, where each $a_i$ is an arbitrary element of $N_{2i}$ and each $b_i$ an arbitrary element of $N_{D(G)-2i}$. If $D(G)$ is odd, then $d_1$ or $d_2$ is even, so $\alpha_j \leq \alpha(G)-\frac{d_1}{2}-\frac{d_2}{2}+\frac{1}{2}=\alpha(G) -\frac{D(G)-3}{2}$. If $D(G)$ is even, then $\alpha_j \leq \alpha(G)-\frac{d_1}{2}-\frac{d_2}{2}+1=\alpha(G) -\frac{D(G)-4}{2}$. In both cases, we have: $\alpha_j \leq \alpha(G) - \left \lfloor \frac{D(G)-3}{2} \right \rfloor$. \\

We say that the cops \emph{occupy} a maximal independent set $S$ of a subgraph $H$ of $G$ if each vertex of $S$ has at least $1$ cop on it. We will show that whenever $\alpha_j$ cops occupy a maximal independent set of $B_i$, we can move them to occupy a maximal independent set of $B_{i+1}$ while keeping $N_{i+2}$ dominated throughout the whole process. Proving this claim obviously finishes the proof of the proposition, since starting by putting the cops on a maximal independent set of $B_0$ and iterating the claim above will force the robber to stay in the subgraph $N_{\geq i+3}$ at the $i$-th iteration, which will lead to his capture. Suppose the cops are occupying a maximal independent set $S$ of $B_i$, then we can successively take every cop that is not in $S \cap (N_{i+1} \cup N_{i+2})$ and move it to a vertex in $N_{i+3}$ that has no edge with any cop located in $N_{i+1} \cup N_{i+2}$. By definition of $\alpha_j$, the cops will eventually occupy a maximal independent set of $B_{i+1}$, and since we are not moving any cop from $S \cap (N_{i+1} \cup N_{i+2})$, the set $N_{i+2}$ stays dominated by the cops through the whole process, which finishes the proof.
\end{proof}

\section{Cop number and domination number}

We now prove \Cref{GammaMinus1} and \Cref{prop:diam-domi} that express relations between cop number and domination number, using the same techniques as their analogues for independence number,  \Cref{AlphaMinus1} and \Cref{prop:diam-idpd}.

\begin{proof}[Proof of \Cref{GammaMinus1}.]
Let $\Gamma$ be a minimum dominating set of $G$. Let $u$, $v$ be vertices such that $\mathrm{dist}(u,v)=D(G)$, $N_i=N_i(u)$ and $N_{\geq i}=N_{\geq i}(u)$. We have $| \Gamma \cap (\{ u \} \cup N_1) | \geq 1$ since $u$ must be dominated, and take one such vertex $w \in \Gamma \cap (\{ u \} \cup N_1)$. We place one cop on each element of $\Gamma_{w} = \Gamma \backslash \{ w \}$. Since $\Gamma_{w}$ dominates $N_{\geq 3}$, the robber must start in $\{ u \} \cup N_1 \cup N_2$. But $N_6$ must be dominated by $\Gamma_{w}$, hence $\Gamma_{w} \cap N_{\geq 5}$ is not empty. Let $c$ be a cop that started in $N_{\geq 5}$. We move the cop $c$ to $w$ and leave the other cops where they are: the robber cannot escape $\{ u \} \cup N_1 \cup N_2$ since $N_3$ is dominated by the cops placed on the vertices of $\Gamma_{w} \cap (N_2 \cup N_3 \cup N_4)$. Hence the robber will be caught and we used $\gamma(G) - 1$ cops, proving $c(G) \leq \gamma(G)-1$ as wanted.
\end{proof}

\begin{figure}[htbp]\centering
    			\includegraphics[height=3cm]{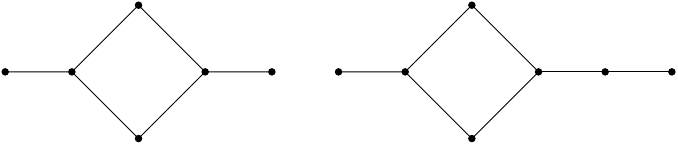}
    			\caption{Examples of connected graphs with $c(G)=\gamma(G)$ for $D(G) \in \{4,5\}$.}
\label{c=g}
	\end{figure}

Like \Cref{AlphaMinus1}, \Cref{GammaMinus1} is optimal in the sense that there are examples of connected graphs with $c(G)=\gamma(G)$ for diameter $D(G) \leq 5$, see \Cref{c=g} for $D(G) \in \{4,5\}$. 

Again, the proof of \Cref{prop:diam-domi} uses a similar idea.

\begin{proof}[Proof of \Cref{prop:diam-domi}.]
Let $k$ be the maximum positive integer such that $m=\left \lfloor \frac{D(G)}{3k+3} \right \rfloor$ satisfies $m-2 \geq k$. Note that $k=\frac{\sqrt{D(G)}}{\sqrt{3}}+O(1)$. Further, let $u$, $v$ be vertices such that $\mathrm{dist}(u,v)=D(G)$. We split the sets $N_i=N_i(u)$ with $0 \leq i \leq (m-1)(3k+3)+3k+2$ into  blocks of size $3k+3$:
\begin{align*}
B_1 &= N_0 \cup N_1 \cup \dots \cup N_{3k+2},\\
B_2 &= N_{3k+3} \cup N_{3k+4} \cup \dots \cup N_{6k+5}, \\
\vdots \\
B_m &= N_{(m-1)(3k+3)} \cup \dots \cup N_{(m-1)(3k+3)+3k+2}.
\end{align*}
We refer to sets of the type $N_{i(3k+3)-1}$, $N_{i(3k+3)}$ and $N_{i(3k+3)+1}$  as \emph{barriers}. For every $i$, let $B'_{i}$ be the block $i$ without the barriers, that is, $B'_{i}=B_{i} \backslash (N_{i(3k+3)} \cup N_{i(3k+3)+1} \cup N_{(i+1)(3k+3)-1})$. Let $\Gamma$ be a dominating set of $G$ of minimal size. For every $j$, $\Gamma \cap (N_j \cup N_{j+1} \cup N_{j+2})$ is non-empty, and hence, for every $i$, $|B'_{i} \cap \Gamma| \geq k$. Let now $S_i \subset B'_{i} \cap \Gamma$ be such that $|S_i| = k$. We put a cop on every vertex of $\Gamma \backslash (S_1 \cup \dots \cup S_m)$. The starting position of the robber must be in one of the $B_i$. On every triplets of consecutive barriers, there is at least one cop. We can then use the at least $m-2$ cops located on triplet of consecutive barriers not adjacent to $B_i$ and put them on each vertex of $S_i$, which is possible since $m-2 \geq k$. Note that while doing that, the robber must stay inside $B_i$ because the barrier $N_{i(3k+3)}$ between $B_{i-1}$ and $B_i$ is still dominated by cops, and so is the barrier $N_{(i+1)(3k+3)}$ between $B_{i}$ and $B_{i+1}$. Consequently, he will be captured by the cops, whose number is

\begin{equation}\label{eq1}
|\Gamma \backslash (S_1 \cup \dots \cup S_m)|= \gamma(G) - km \leq \gamma(G) -k^2 - 2k.    
\end{equation} By inserting $k=\frac{\sqrt{D(G)}}{\sqrt{3}}+O(1)$ into \eqref{eq1}, we get: 
$$c(G) \leq \gamma(G)-\dfrac{D(G)}{3} +O(\sqrt{D(G)}),$$
as desired. 
\end{proof}

\section{Concluding remarks and open problems}

Even if \Cref{AlphaMinus1} is optimal in the sense that there exists graphs of small diameter satisfying $c(G)=\alpha(G)$, we had not found any examples of connected graphs satisfying this equality with $\alpha(G) > 2$. We are wondering whether the following holds:

\begin{Question} \label{Q1}
Given a positive integer $n$, does there exist a connected graph $G$ such that $c(G)=\alpha(G) \geq n$?
\end{Question}

The problem is the same with \Cref{GammaMinus1} and graphs satisfying $c(G)=\gamma(G)$. We had not found any examples of connected graphs satisfying this equality with $\gamma(G) > 3$:

\begin{Question} \label{Q2}
Given a positive integer $n$, does there exist a connected graph $G$ such that $c(G)=\gamma(G) \geq n$?
\end{Question}

Anurag Bishnoi and Jérémie Turcotte independently brought to our attention examples of connected graphs satisfying the equalities in \Cref{Q1} and \Cref{Q2} for more small values of $n$. For instance, the Paley graph on $17$ vertices satisfies $c(G)=\alpha(G)=3$. They also provided us with examples of connected graphs satisfying $c(G)=\gamma(G) = n \in \{4,5,6,7\}$. For example, for $n=7$, one may consider the Hoffman Singleton graph, whose cop number can be easily derived using a classic result due to Aigner and Fromme \cite{AIGNER} that the cop number of a graph with girth at least $5$ is at least its minimum degree.

The questions still remain open for $\alpha \geq 4$ and $\gamma \geq 8$.

\section*{Acknowledgement}

The authors would like to thank Professor Béla Bollobás for his valuable comments.

\bibliographystyle{abbrvnat}  
\renewcommand{\bibname}{Bibliography}
\bibliography{bibliography}

\begin{thebibliography}{4}
\providecommand{\natexlab}[1]{#1}
\providecommand{\url}[1]{\texttt{#1}}
\expandafter\ifx\csname urlstyle\endcsname\relax
  \providecommand{\doi}[1]{doi: #1}\else
  \providecommand{\doi}{doi: \begingroup \urlstyle{rm}\Url}\fi

\bibitem[Aigner and Fromme(1984)]{AIGNER}
M.~Aigner and M.~Fromme.
\newblock A game of cops and robbers.
\newblock \emph{Discrete Applied Mathematics}, 8\penalty0 (1):\penalty0 1--12,
  1984.
\newblock ISSN 0166-218X.
\newblock \doi{https://doi.org/10.1016/0166-218X(84)90073-8}.
\newblock URL
  \url{https://www.sciencedirect.com/science/article/pii/0166218X84900738}.

\bibitem[Nowakowski and Winkler(1983)]{Now}
R.~J. Nowakowski and P.~Winkler.
\newblock Vertex-to-vertex pursuit in a graph.
\newblock \emph{Discret. Math.}, 43\penalty0 (2-3):\penalty0 235--239, 1983.
\newblock \doi{10.1016/0012-365X(83)90160-7}.
\newblock URL \url{https://doi.org/10.1016/0012-365X(83)90160-7}.

\bibitem[Quillot()]{Qui}
A.~Quillot.
\newblock Jeux et pointes fixes sur les graphes.
\newblock \emph{Ph.D. Dissertation, Université de Paris VI, 1978}.

\bibitem[Turcotte(2022)]{2k2free}
J.~Turcotte.
\newblock Cops and robbers on 2{K}2-free graphs.
\newblock \emph{Discrete Mathematics}, 345\penalty0 (1):\penalty0 112660, 2022.
\newblock ISSN 0012-365X.
\newblock \doi{https://doi.org/10.1016/j.disc.2021.112660}.
\newblock URL
  \url{https://www.sciencedirect.com/science/article/pii/S0012365X21003733}.

\end{thebibliography}

\end{document}